\documentclass{article}

\usepackage{amsmath,amssymb,booktabs,enumitem,graphicx,tikz}
\usepackage{natbib}
\usetikzlibrary{arrows.meta,positioning}
\title{Decision-Making in Digital Twins: Optimization, Control, and Learning}
\author{Enlu Zhou \\
H. Milton Stewart School of Industrial and Systems Engineering \\
Georgia Institute of Technology}
\date{}

\begin{document}
\maketitle

\begin{abstract}
A digital twin is a virtual representation of a physical system, continuously updated with data from the physical system, enabling dynamic information exchange and decision-making. This bidirectional interaction, which distinguishes digital twins from conventional simulation models, introduces new challenges to simulation-based optimization, control, and learning. 
This paper reviews recent developments in decision-making methodologies motivated by digital twins. We focus on three interconnected areas: simulation optimization, sequential decision making, and risk-aware decision making. For simulation optimization, we discuss data-driven simulation optimization under input uncertainty, streaming input data, decision-dependent observations, and limited computational and data-collection resources. For sequential decision making, we review the connections among digital twins, stochastic optimal control, and reinforcement learning, with emphasis on Bayesian approaches that update models and policies as new information becomes available. For risk-aware decision making, we discuss distributionally robust and Bayesian risk formulations as alternative approaches for addressing model uncertainty. Throughout the paper, we emphasize the common problem of making decisions while the underlying model is being updated. We conclude by identifying research  challenges and opportunities.
\end{abstract}


\section{Introduction}
A digital twin is ``a set of virtual information constructs that mimics the structure, context, and behavior of a natural, engineered, or social system (or system of systems), is dynamically updated with data from its physical twin, has a predictive capability, and informs decisions that realize value'' \cite{nasem2024digital}. Digital twins provide an appealing tool for modeling, comprehending, and optimizing complex real-life systems. They have significant potential in manufacturing, transportation, supply chains, energy systems, healthcare, and other engineering and scientific domains \cite{kritzinger2018digital,jones2020characterising,heluany2023survey,wu2023comprehensive}.

At the core of many digital twins is a stochastic simulation model or a system of such models. The defining feature of a digital twin, distinguishing it from conventional simulation, is its bidirectional interaction with the physical system. The physical system continually provides its digital twin with real-time data, enhancing the fidelity of the simulation models, while predictions and decisions based on these models are provided back to the physical system. This interaction produces a closed-loop process involving data, models, and decisions \cite{biller2022simulation,matta2023digital,zhou2024tutorial}.

Conventional simulation-based decision making usually assumes access to a fixed simulation model. The model may be estimated from historical data, but it is typically treated as unchanged during the decision-making process. Digital twins challenge this assumption. New observations arrive over time, the estimated input distributions may change, decisions may affect future observations, and decisions often must be updated before the underlying optimization problem can be solved to high accuracy. The digital twin is therefore not only a simulation model. It is an adaptive decision-making system.

Most broad reviews of digital twins focus on definitions, architectures, enabling technologies, communication mechanisms, or application domains \cite{kritzinger2018digital,jones2020characterising,heluany2023survey,wu2023comprehensive}. These reviews provide valuable descriptions of how digital twins are constructed and used, but decision-making methodologies have received comparatively less unified treatment. In particular, simulation optimization,  stochastic control, and reinforcement learning are often studied by different research communities, even though they address closely related questions within digital twins.

This review is organized around the following central question:
\begin{quote}
{\it How should a digital twin make and update decisions while learning from continuously arriving data?}
\end{quote}
We study this question through three areas. The first is \emph{data-driven optimization}, including simulation optimization with input uncertainty, streaming observations, and nonstationarity. The second is \emph{sequential decision making}, including stochastic optimal control and reinforcement learning. The third is \emph{risk-aware decision making}, which cuts across both optimization and sequential decision problems.

The scope of the paper is focused on methodological developments that connect data to decisions. This focused scope allows us to review a manageable body of literature while developing a unified perspective on simulation optimization, stochastic control, and reinforcement learning, with particular attention to risks arising from the sim-to-real gap.

\section{Digital Twins as Adaptive Decision-Making Systems}
A useful abstraction of a digital twin consists of five interacting components: (i) a physical system; (ii) a stream of observations from that system; (iii) one or more models that constitute the twin; (iv) a decision-making module; and (v) actions implemented in the physical system. Together, these components form a closed loop: observations update the models, the models predict the consequences of candidate actions, and the decision-making module selects an action to implement. That action changes the physical system and may also change which observations arrive next. Data, models, and decisions can therefore become endogenously coupled.

Within this loop, simulation is the operational core that turns the twin's representation into predictions and decisions. \cite{biller2022simulation} describe simulation as a critical technology in digital-twin development, while application-oriented frameworks, such as the O2DES.NET digital twins studied by \cite{li2020o2des}, illustrate how simulation supports operational decisions in ports and warehouses. These works explain how a twin can represent a physical system. 
The distinctive methodological challenge arises because that representation is repeatedly updated as new data arrive.

To make this dynamic explicit, let $M_t$ denote the physical system's current environment and $\widehat{M}_t$ denote the digital twin's current model after assimilating the data available at synchronization stage $t$. Conventional simulation-based optimization or decision-making often fits an input distribution from a fixed batch of historical data, runs one model, and obtains a decision. A digital twin instead uses each new batch of observations generated from $M_t$ to update $\widehat{M}_t$ to $\widehat{M}_{t+1}$, conducts additional simulations under the updated model, and updates the decision. Simulation outputs generated at different stages may therefore follow different distributions. They may also be correlated because successive model estimates are based on overlapping data streams \cite{zhou2024tutorial}.

The evolving sequence $\{\widehat{M}_t\}$ turns the twin from a collection of simulation models into an adaptive decision process. Let $S_t$ denote the operational state of the physical system at stage $t$. The twin uses $\widehat{M}_t$ and $S_t$ to select a decision $x_t$, implements $x_t$ in the physical system, and receives observations from $M_t$ that inform the next model update. This coupling raises three general questions: how should information from different model versions be combined; how should limited data-collection and simulation resources be allocated; and how should decisions be updated when the model evolves? These questions motivate the data-driven optimization methods in Section~3 and the sequential decision-making methods in Section~4. They also motivate the risk-aware formulations in Section~5, since each action is selected using an imperfect and evolving representation of the physical system.

\begin{figure}[!ht]
\centering
\resizebox{0.96\linewidth}{!}{%
\begin{tikzpicture}[
  font=\small,
  >=Latex,
  stage/.style={draw=blue!55!black, fill=blue!6, rounded corners=2pt,
    align=center, minimum height=1.15cm, text width=2.65cm, line width=0.65pt},
  decide/.style={draw=orange!70!black, fill=orange!9, rounded corners=2pt,
    align=center, minimum height=1.15cm, text width=2.65cm, line width=0.65pt},
  flow/.style={->, line width=0.85pt, draw=teal!70!black},
  history/.style={draw=gray!65, fill=gray!7, rounded corners=2pt,
    align=center, text width=5.2cm, font=\footnotesize}
]
\node[stage] (physical) {\textbf{1. Physical system} $M_t$\\operational state $S_t$};
\node[stage, right=1.0cm of physical] (data) {\textbf{2. Observe}\\stream data and estimate inputs};
\node[stage, right=1.0cm of data] (twin) {\textbf{3. Update twin}\\model $\widehat{M}_t$};
\node[decide, right=1.0cm of twin] (decision) {\textbf{4. Decide}\\simulate, predict, and select $x_t$};

\draw[flow] (physical) -- (data);
\draw[flow] (data) -- (twin);
\draw[flow] (twin) -- (decision);
\draw[flow] (decision.north) -- ++(0,1.05cm) coordinate (feedbacktop) -- node[midway, above=0.14cm, font=\normalsize] {implement $x_t$; system evolves to $t+1$} (physical.north |- feedbacktop) -- (physical.north);

\node[history, below=1.2cm of data, xshift=1.8cm] (history)
  {Historical observations and simulation outputs\\may be retained and reused across model updates};
\draw[->, dashed, draw=gray!75, line width=0.65pt] (history.north west) to[bend left=10] (data.south);
\draw[->, dashed, draw=gray!75, line width=0.65pt] (twin.south) to[bend left=10] (history.north east);
\node[text=gray!75!black, below=0.08cm of history]
  {overlapping data can induce cross-stage dependence};
\end{tikzpicture}%
}
\caption{Dynamic decision making with a digital twin. At each stage, physical observations update the input model and simulator; the updated twin generates predictions for selecting an action; and implementation changes the physical system and its future data. Reusing observations across stages can make successive model estimates and simulation outputs dependent.}
\label{fig:dt-decision-loop}
\end{figure}
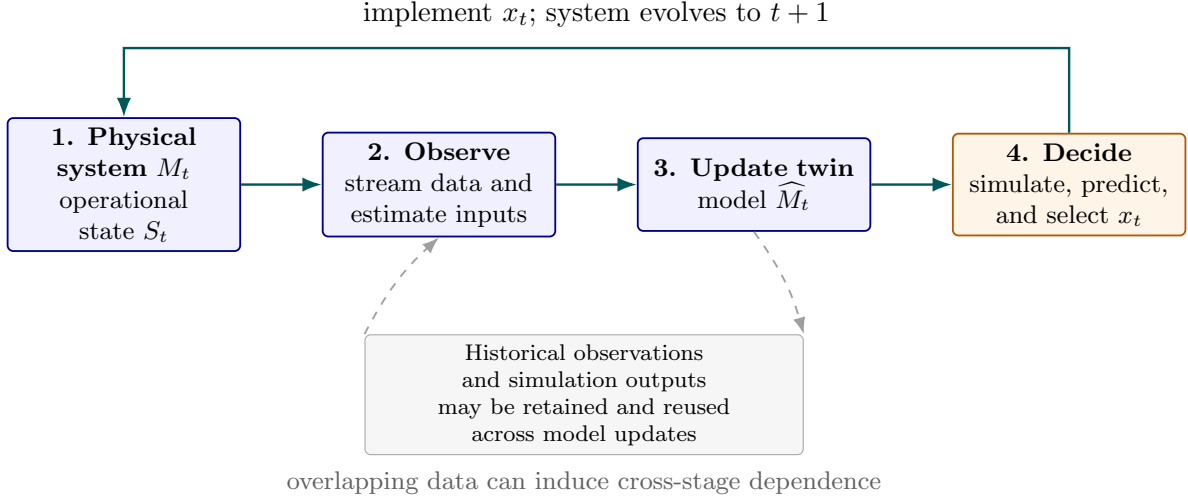

Figure~\ref{fig:dt-decision-loop} summarizes this data--model--decision loop and its repeated synchronization with the physical system.

\section{Data-Driven Simulation Optimization in Digital Twins}
Section~2 framed a digital twin as a recurring data--model--decision loop in which physical observations update $\widehat{M}_t$ and the updated twin selects $x_t$. This section focuses on the optimization component of that loop: we first isolate input-model uncertainty at a synchronization stage, then consider streaming model updates and the joint allocation of data-collection and simulation effort, and finally discuss context-dependent decisions for rapid online deployment.

\subsection{Simulation optimization with input uncertainty}
Section~2 represents the physical system and its digital twin by $M_t$ and $\widehat{M}_t$, respectively. We first focus on one important source of discrepancy between them: input-model error. Specifically, suppose that the simulator's response mapping and all other model components are correctly specified, so that the only difference between $M_t$ and $\widehat{M}_t$ is the distribution of the stochastic inputs. Let $F_t$ denote the true input distribution governing $M_t$, and let $\widehat{F}_t$ denote the data-driven input distribution used by $\widehat{M}_t$. Under this restriction, we may identify $M_t$ with $F_t$ and $\widehat{M}_t$ with $\widehat{F}_t$, thereby isolating the part of the twin that is estimated and updated from physical data. When analyzing a single synchronization stage, we suppress the time subscript and write $F$, $\widehat{F}$, and $x$ for $F_t$, $\widehat{F}_t$, and $x_t$, respectively. The operational decision problem under the physical system is then
\begin{equation}
\min_{x\in\mathcal{X}} \; \mathbb{E}_{F}[H(x,\xi)],
\label{eq:so}
\end{equation}
where $H$ is the correctly specified response mapping of the simulation model, and $\xi\sim F$ is the physical system's stochastic input. Since $H$ generally does not have an explicit form and can only be evaluated by running simulations, solving problem~\eqref{eq:so} is often called simulation optimization \cite{fu2015handbook}. The twin does not observe $F$ directly; it simulates under $\widehat{F}$, which is estimated from the available physical data. Consequently, even when every other part of the simulator is correct, a discrepancy between $\widehat{F}$ and $F$ can cause a decision optimized in $\widehat{M}_t$ to perform poorly when implemented in $M_t$.

A particularly important special case of problem~\eqref{eq:so} is ranking and selection (R\&S) \cite{hong2021review}. In R\&S, the decision space $\mathcal{X}$ is a finite set of categorical alternatives, often called systems or designs, and simulation is used to estimate and compare their expected performances in order to select the best one. Thus, unlike continuous simulation optimization, the decision is which alternative to select rather than the value of a continuous decision vector.

Considering the input model error, two sources of uncertainty exist in digital twins. Simulation uncertainty is caused by stochastic simulation outputs and limits how precisely candidate decisions can be compared under $\widehat{F}$; it can generally be reduced by increasing the number of replications. Input uncertainty arises because $\widehat{F}$ is estimated from limited observations generated by $M_t$ and may differ from $F$. Additional simulation replications reduce simulation uncertainty but do not close this input-model gap, which instead requires more informative physical data or improved input modeling. Distinguishing the two sources is therefore essential when assessing confidence in a decision produced by the digital twin.

The input-uncertainty literature therefore provides a statistical foundation for calibrating and validating digital twins. Methods such as Bayesian, bootstrap, subsampling, and nested-simulation quantify how the uncertainty in  finite input data propagates to simulation outputs and construct confidence intervals or risk measures for performance under input uncertainty \cite{chick2001input,zouaoui2003accounting,zouaoui2004accounting,xie2014bayesian,barton2014quantifying,lam2022subsampling,zhu2020risk}. Their common message for digital twins is that simulation noise and input-estimation error must be separated when assessing decision reliability and deciding whether more physical data or more simulation replications are needed.

Input uncertainty describes the discrepancy that the decision maker must confront, but it does not by itself prescribe an attitude toward that uncertainty. To hedge againt risk associated with input uncertainty, two classes of methods exist: Distributionally robust optimization (DRO), and Bayesian risk optimization (BRO). DRO replaces the estimated input distribution by an ambiguity set and optimizes performance under its worst-case input distribution \citep{gao2017robust,xiao2018simulation,xiao2020optimal,fan2020distributionally}. This provides protection against model error, although it can be overly conservative depending on how the ambiguity set and its radius are chosen. BRO, by contrast, places a posterior distribution on the unknown input model and applies a risk functional to the resulting posterior distribution of performance \citep{wu2018bayesian}. Posterior expectation, quantiles, and tail-risk measures span attitudes from average-case to increasingly conservative behavior, while Bayesian updating assimilates new data. Thus, BRO separates learning about the input model from choosing an attitude toward the remaining uncertainty. Section~5.1 reviews these criteria and their relationship. 


\subsection{Simulation optimization with streaming input data and resource allocation}
While the previous subsection considers a fixed synchronization stage, 
a critical aspect of digital twins is the synchronization of the physical system and its digital twin through streaming data. Streaming observations arrive sequentially, possibly in random batches of varying sizes. At each stage, observations generated by $M_t$ update the estimated input distribution within $\widehat{M}_t$. Because simulation time is limited, only a finite number of new replications can be generated under each updated model.

\cite{wu2024datadriven} considered a fixed-confidence R\&S problem with streaming input data. Their objective is to identify the true best design with a prescribed probability of correct selection. Simulation outputs generated at different stages are neither independent nor identically distributed because the estimated input parameters change over time and are themselves correlated. To address this issue, they proposed a moving-average estimator that retains only a recent portion of the accumulated simulation outputs. The amount of retained data controls a trade-off: older outputs help reduce simulation uncertainty but may have larger bias because they were generated under less accurate input models. The resulting sequential-elimination procedures provide decisions with either an exact or asymptotically valid guarantee of correct selection.

\cite{wang2022fixed,wang2024optimal} studied a related fixed-budget setting. Rather than minimizing simulation effort subject to a probability-of-correct-selection requirement, they allocate a fixed stage-wise simulation budget to maximize the probability of correct selection. This formulation is useful when a digital twin must generate an updated decision within a prescribed time window.

\cite{he2024stochastic} considered multi-period continuous simulation optimization with streaming input data. At each period, newly collected data update the simulator's input-parameter estimates, and the decision must be updated before the next batch arrives. They developed a multi-period stochastic-approximation framework with restart and warm-start variants. The warm-start method takes fewer and shorter gradient steps as the input estimates become more precise, achieving the optimal convergence rate while substantially reducing computation across periods.

\cite{liu2021bayesian,liu2024bayesian} considered a similar setting of continuous stochastic optimization with streaming input data. Their Bayesian stochastic gradient descent method jointly updates the posterior distribution of an unknown input parameter and the decision variable. For decision-independent data, they showed that the method asymptotically achieves convergence rates comparable to classical stochastic gradient descent. More importantly for digital twins, they also considered decision-dependent observations. In this setting, data generated at time $t$ depend on the current decision $x_t$, so the observations are not identically distributed. They established posterior consistency under suitable distinguishability conditions and convergence of the decision sequence to an optimal solution. This work provides a direct connection between online estimation and online optimization in an evolving simulation model.

The preceding studies treat the incoming data stream and the stage-wise simulation budget as given. In a digital twin, however, data acquisition and simulation may compete for limited resources or proceed concurrently. Waiting until sufficient input data have been collected before initiating simulation can unnecessarily delay decisions, whereas allocating effort to simulation too early may rely on an inaccurate input model. This leads to a joint resource-allocation problem: how should effort be divided between learning the input model and reducing simulation uncertainty? The literature can be organized by whether data collection and simulation share a common budget or separate resources, and whether the two activities occur sequentially or concurrently.

\cite{wu2017ranking} formulated budget allocation under input uncertainty as a two-stage problem. \cite{xu2020joint} studied joint allocation of resources to input-data collection and simulation while accounting for their respective costs. \cite{kim2022optimizing} proposed a sequential acquisition procedure based on the posterior probability that a system is best.

\cite{wang2025ranking} considered simultaneous allocation when data collection and simulation use separate resources and can proceed in parallel. Their allocation policies determine how many observations should be collected from each input source and how many simulation replications should be assigned to each design. The interaction is important: additional input data reduce model uncertainty, whereas additional simulation replications reduce simulation uncertainty. Large-deviations analysis yields balance conditions for allocating both resources, and the resulting sequential algorithms seek to achieve these conditions adaptively. This formulation is natural for digital twins, where data collection and simulation often occur concurrently.

\subsection{Context-dependent decisions and online applications}
Not all data are used to estimate input distributions. Some data represent observable contexts, covariates, or side information. In this case, the decision maker seeks a policy that maps the observed context to a decision rather than a single robust decision.

A closely related literature studies contextual stochastic optimization, in which historical context--outcome pairs are used to learn a decision rule for a downstream stochastic program; see \cite{sadana2025survey} for a recent overview. That literature generally exploits analytical structure of the objective function that makes the optimization problem computationally tractable. Context-dependent simulation optimization addresses a different regime: a candidate decision's performance is available only through noisy and potentially expensive black-box simulation, so a tractable optimization oracle generally does not exist. This distinction delineates the contextual simulation optimization methods reviewed below.

\cite{shen2021ranking} studied ranking and selection with covariates for personalized decisions. \cite{li2020context,li2022efficient} developed Bayesian methods for learning and optimizing context-dependent designs. \cite{cakmak2024contextual} combined Gaussian-process modeling with optimal computing-budget allocation for contextual ranking and selection, by exploiting correlations across contexts and alternatives. 


More recent work broadens both the decision objective and the structure used to share information across contexts. \cite{du2024contextual} derive rate-optimal allocation rules for small and large context spaces in personalized medicine, concentrating simulation on context--treatment pairs that are hardest to distinguish. \cite{keslin2024ranking} propose ranking and contextual selection, which chooses an offline experimental design and uses classification for fast online decisions at newly observed covariates, together with pointwise probability-of-good-selection and postexperiment optimality-gap guarantees. \cite{zhang2025topm} extend contextual selection from a single best design to the top-$m$ designs for every context and develop an asymptotically optimal sequential sampling policy. Moving from finite or unstructured contexts to continuous structure, \cite{li2026quadratic} exploit quadratic response models to preselect informative context points and allocate simulation efficiently over a continuous context space. \cite{hong2019offline} and \cite{jiang2020online} developed an offline-simulation-online-application paradigm, where a large number of scenarios are simulated offline and the resulting data are used to train a predictive rule; when a new context is observed online, the rule provides a decision without requiring additional expensive simulation.

These developments make context-dependent simulation optimization especially relevant to digital twins: substantial computation can be performed offline, while online contextual observations are converted quickly into operational decisions that fits the current context. 

\section{Sequential Decision Making in Digital Twins}
\subsection{From one-time optimization to sequential decisions}
Section~3 largely considers a decision optimized for the digital twin's current model. In many operational settings, however, the current decision changes the physical state, future costs, and the observations available for later model updates. Decision quality must therefore be assessed over a horizon rather than through immediate performance alone. A Markov decision process (MDP) formalizes this dependence: the physical system constitutes the environment, its current condition is summarized by a state, and a policy maps that state to an action whose consequences unfold over time. The MDP can be viewed as a decision abstraction supported by the digital twin.

To connect with the notations introduced in Section~2, let $S_t$ denote the physical system's operational state at synchronization stage $t$. The true environment is $M_t$, whereas $\widehat{M}_t$ is its digital-twin approximation. A policy $\pi_t$ maps the state to the implemented action, $A_t=\pi_t(S_t)$; hence, the generic decision $x_t$ in Section~2 is the exercised action $A_t$ at stage $t$. In the digital twin, given a state--action pair $(s,a)$, $\widehat{M}_t$ can generate the next state, while data and state estimation produce the current state, and physical implementation produces the next observation and realized performance. A one-stage reward or cost may be embedded in the digital twin, or it may be specified externally from engineering objectives, economic criteria, safety penalties, and operating constraints, and then evaluated on simulated trajectories. Thus, the digital twin simulates the transition from a state-action pair to the next state and provides the information needed to evaluate candidate actions and policies. This correspondence also depends on synchronization frequency and observability. When synchronization and decision epochs coincide, $\widehat{M}_t$ is updated before each action; when synchronization is less frequent, the same twin model is held fixed across the intervening MDP steps. If the physical state is not observed directly, sensor observations and a state estimator instead define an information or belief state, and a partially observable MDP (POMDP) may be more appropriate. 

Specifically, let $P_{\theta_t}$ denote the transition kernel simulated by the digital twin and parametrized by $\theta_t$. A policy $\pi$ is defined as a sequence $\{\pi_0, \pi_1, \ldots\}$, where each $\pi_t$ maps the state $S_t$ either to a single action $A_t$ (a deterministic policy) or to a probability distribution over the feasible actions (a randomized policy). Then, the controlled sate evolution in the digital twin can be written as
\begin{equation}
S_{t+1} \sim P_{\theta_t}(\cdot \mid S_t,A_t),
\qquad
A_t \sim \pi_t(\cdot \mid S_t).
\label{eq:mdp-dynamics}
\end{equation}
The parameter $\theta_t$, or a posterior distribution over it, is updated as physical observations arrive.  For a fixed model $\theta$, the finite-horizon $\pi$ can be evaluated through
\begin{equation}
V_{\theta}^{\pi}(s)
=
\mathbb{E}_{\theta}^{\pi}
\left[
\sum_{t=0}^{T-1} \gamma^t c_t(S_t,A_t)
+ \gamma^T c_T(S_T)
\,\middle|\, S_0=s
\right],
\label{eq:policy-objective}
\end{equation}
where $c_t$ is the stage-wise cost and $\gamma\in (0,1)$ is a discount factor. For an infinite-horizon discounted MDP with stationary cost functions, the corresponding criterion is $V_{\theta}^{\pi}(s)=\mathbb{E}_{\theta}^{\pi}[\sum_{t=0}^{\infty}\gamma^t c(S_t,A_t)\mid S_0=s]$. For such an infinite-horizon MDP, under standard regularity conditions, such as finite state and action spaces and bounded costs, there exists an optimal stationary Markov policy whose action rule does not depend explicitly on time.

MDPs provide a common mathematical framework for the two main approaches developed in the following subsections. Stochastic optimal control typically begins with an explicit model and computes policies by dynamic programming or related methods, whereas reinforcement learning seeks good policies through interaction with an environment \citep{bertsekas2017dynamic,sutton2018reinforcement}. A digital twin can support both routes: it provides an evolving predictive model for control and a virtual environment for policy learning, while new physical observations update the model used by either approach. For the infinite-horizon discounted problem with known $\theta$, the value of an optimal stationary policy satisfies the Bellman equation
\begin{equation}
V_{\theta}^{\star}(s)
=
\min_{a\in\mathcal{A}(s)}
\left\{
c(s,a)
+ \gamma \mathbb{E}_{S'\sim P_{\theta}(\cdot\mid s,a)}
\left[V_{\theta}^{\star}(S')\right]
\right\}.
\label{eq:bellman}
\end{equation}
Equations~\eqref{eq:mdp-dynamics}--\eqref{eq:bellman} expose the central challenge for digital twins: the transition model used for evaluation and control  is no longer fixed as $P_\theta$ in the digital twin; instead, itself is  updated and evolving as $P_{\theta_t}$ while decisions are being made.

This evolving-model challenge can be addressed through explicit model control, direct policy learning, or a hybrid of the two. As an illustration of the hybrid route, \cite{schena2024reinforcement} proposed ``reinforcement twinning,'' which combines data assimilation and model-based optimal control with model-free reinforcement learning in a virtual environment before policies interact with the physical system. The next two subsections develop the two constituent routes in turn: Section~4.2 reviews learning-based and adaptive stochastic control, and Section~4.3 reviews simulation-based reinforcement learning.

\subsection{Learning-based and adaptive stochastic control}
Stochastic control provides a direct solution to sequential decision making when the digital twin contains an explicit dynamics model. Rather than learning a policy only from trial-and-error interaction, a controller repeatedly solves a finite- or infinite-horizon optimization problem and applies feedback to the physical system. The central challenge is that the model, disturbance distribution, and sometimes the control objective are only partially known and continue to change as observations arrive.

One family of methods learns a probabilistic dynamics model and propagates its uncertainty through planning. \cite{deisenroth2011pilco} developed PILCO, which uses a Gaussian-process dynamics model and incorporates predictive uncertainty into long-term policy evaluation. Its data efficiency illustrates a key advantage of explicit model learning for digital twins, where physical experiments may be costly but the learned model can be simulated repeatedly. At the same time, planning against an estimated model can bias the controller toward regions in which the twin is inaccurate, so uncertainty representation must be coupled to the control calculation. 

Model predictive control provides a natural architecture for this coupling. At each decision time, MPC solves a constrained finite-horizon problem using the current model, applies the first action, and replans after receiving new measurements. The survey of \cite{hewing2020learning} organizes learning-based MPC around learning the dynamics, learning controller components such as costs and constraints, and using MPC to provide safety guarantees for learning-based policies. As a concrete example, \cite{koller2018safe} constructed Gaussian-process confidence sets for predicted trajectories and combined them with terminal-set constraints to guarantee high-probability safety during exploration. 

Bayesian methods make the learning-control interaction explicit through a posterior over models. \cite{wabersich2020bayesian} combined posterior sampling with model predictive control, using sampled dynamics models to explore efficiently while retaining the computational structure of MPC. \cite{li2026bayesianmpc} developed a Bayesian risk-averse MPC framework with consistency and stability guarantees, explicitly incorporating posterior model uncertainty and risk sensitivity into receding-horizon control. A related episodic formulation is especially appropriate when a digital twin is synchronized periodically rather than continuously. \cite{shapiro2025episodic} update the posterior distribution of the unknown exogenous randomness at the beginning of each episode, solve a Bayesian-average control problem, implement the resulting policy, and use the new observations in the next episode. For convex costs and linear dynamics, the method can be implemented with stochastic dual dynamic programming.

Taken together, these approaches define a broad control-side role for digital twins: update an explicit predictive model, propagate its uncertainty into planning, enforce operational constraints during adaptation, and choose an update frequency that matches the physical system and available computation. When repeatedly solving an explicit control problem is impractical, the twin can instead support direct policy learning.

\subsection{Simulation-based reinforcement learning}
Whereas the preceding approaches use an explicit dynamics model inside an optimization problem, simulation-based reinforcement learning uses the digital twin as an interactive training environment. Policies can interact with the simulator far more cheaply, quickly, and safely than with the physical system. Simulation can expose a learning agent to rare failures, explore risky actions without operational risk, and generate possible trajectories in parallel. 

A first approach treats the simulator as a learned or specified transition model and alternates between physical experience and simulated planning. The classical Dyna architecture of \cite{sutton1990integrated} interleaves model learning, policy learning, and planning with synthetic transitions. Modern model-based reinforcement learning follows the same principle at larger scale. For example, \cite{janner2019trust} showed that short model-generated rollouts branched from real states can improve sample efficiency while limiting the bias that accumulates when a learned model is rolled out for too long. This trade-off is especially important for digital twins: inexpensive virtual experience is valuable only to the extent that the twin remains predictive in the states visited by the policy.

A second approach trains primarily in simulation and then transfers the policy to the physical system. The central difficulty is the sim-to-real gap, since a policy may exploit modeling errors that are harmless in the simulator but consequential after deployment. Dynamics randomization addresses this problem by training over a distribution of plausible environments rather than a single calibrated model. \cite{peng2018simtoreal} demonstrated that randomizing physical parameters can produce control policies that transfer to real hardware. \cite{chebotar2019closing} closed the loop by adapting the simulation-parameter distribution using a small number of real-world rollouts, thereby concentrating training on simulations that better reproduce physical behavior.

Digital twins extend these ideas from a one-time transfer problem to a continuing learning loop. For instance, \cite{voogd2023reinforcement} combined virtual and real data across model-, hardware-, and vehicle-in-the-loop stages for autonomous driving, together with domain randomization and adaptation. 
Together, these strands suggest a simulation-based reinforcement-learning cycle for digital twins: learn and test policies in the simulator, deploy cautiously, and use new physical trajectories to update both the twin and the policy. The methods above concern how the twin supplies experience and supports policy learning. When its transition model remains uncertain, an additional question is how that uncertainty should enter policy evaluation and optimization: distributionally robust and Bayesian risk-aware formulations will be examined in Sections~5.2 and~5.3. Repeated updating also raises whether trajectories generated under earlier policies or earlier twin versions can be reused, which is the focus of the next subsection.

\subsection{Reusing experience across model updates}
Repeatedly updating a digital twin creates both a difficulty and an opportunity. The difficulty is that data were generated under earlier policies and earlier versions of the model. The opportunity is that past trajectories or evaluations may provide valuable information for the current policy evaluation.

A useful distinction is whether historical data are mismatched because the policy changed, the environment model changed, or both. \cite{lin2025reusing} studied the reuse of historical trajectories in natural policy gradient methods through importance sampling, providing convergence and rate guaranties for incorporating past trajectories into current policy updates. When the environment model changes, \cite{tirinzoni2019transfer} developed multiple-importance-sampling estimators for transferring complete trajectories between environments, including safeguards against negative transfer; \cite{chandak2020optimizing} used nonuniform weights on past data to optimize forecasts of future policy performance in nonstationary MDPs; and \cite{liu2023reusing} proposed a regression-assisted doubly robust estimator that remains asymptotically unbiased. When the environment is fixed but the policy evolves, \cite{queeney2021generalized} extended proximal policy optimization to reuse data from several recent policies while retaining an off-policy improvement guarantee. Together, these corrections for model shift and policy shift can reduce the cost of repeatedly resolving sequential decision problems after each digital-twin update.

\section{Risk-Aware Decision Making}
Section~3 introduced input model uncertainty as a discrepancy between the physical system and the digital twin, while Section~4 mainly concerned the discrepancy in transition models for sequential decision making. This motivates us to think about the question: how should the decision criterion respond to uncertainty in those models? Risk may arise from the intrinsic randomness of the physical system or from incomplete knowledge of its model. These two sources are often called aleatoric and epistemic uncertainty, respectively. In a digital twin, epistemic uncertainty changes over time because the twin continually receives data; hence, we focus on epistemic uncertainty in this review. A risk-aware decision-making method should protect against consequential model error while allowing the degree of protection to evolve with the amount and quality of available information. 

\subsection{Distributionally robust optimization and Bayesian risk optimization}
Distributionally robust optimization (DRO) is a principal framework for decisions under distributional ambiguity. Instead of committing to a single estimated distribution, DRO optimizes performance under the worst distribution in an ambiguity set. Let $\mathcal{P}_n$ denote an ambiguity set constructed from the $n$ observations currently available to the twin. The DRO counterpart of problem~\eqref{eq:so} is
\begin{equation}
\min_{x\in\mathcal{X}}
\;\sup_{F\in\mathcal{P}_n}
\mathbb{E}_{F}[H(x,\xi)].
\label{eq:dro}
\end{equation}
Several constructions of $\mathcal{P}_n$ balance statistical justification and tractability. Moment-based sets use confidence regions for moments or conic descriptions of moment and support information \citep{delage2010distributionally,wiesemann2014distributionally}. For sample-based models, Wasserstein balls and optimal-transport duality yield tractable formulations, out-of-sample guarantees, and structural results \citep{esfahani2018data,blanchet2019quantifying,gao2023distributionally}. These tools are attractive for digital twins whose current models are estimated from limited physical data.

The strength of DRO is also the source of its conservatism. Once an ambiguity set is specified, the criterion treats its least favorable distribution as decisive. The geometry and radius of the set therefore encode both statistical uncertainty and risk attitude, although these are conceptually different. A radius calibrated for high-confidence coverage can yield excessive pessimism, while a radius chosen for good nominal performance may provide little protection. Repeated data assimilation also requires the ambiguity set to be reconstructed or recalibrated as the twin evolves. Thus, data-driven DRO adapts to data through the ambiguity set, but its standard worst-case objective can still be overly conservative given the available data.

Bayesian risk optimization (BRO) provides a more flexible alternative to this worst-case construction. Let $\mu_n$ denote the posterior distribution on the unknown input distribution $F$ after $n$ observations, and let $\rho_{\mu_n}$ be a risk functional with respect to that posterior. A generic BRO problem is
\begin{equation}
\min_{x\in\mathcal{X}}
\;\rho_{\mu_n}\!\left(
\mathbb{E}_{F}[H(x,\xi)]
\right).
\label{eq:bro}
\end{equation}
Thus, BRO applies risk to the posterior distribution of performance rather than taking the worst case over a fixed ambiguity set \citep{wu2018bayesian}. Posterior expectation yields Bayesian-average optimization; posterior quantiles and conditional value-at-risk emphasize adverse but posterior-plausible models; and the essential supremum approaches a worst-case criterion. At the same time, Bayes' rule gives a natural mechanism for assimilating each new batch of physical observations. As the posterior concentrates, epistemic uncertainty and the associated risk adjustment contract automatically. \cite{wu2018bayesian} established statistical properties of BRO, \cite{cakmak2021solving} developed nested stochastic gradient descent methods for solving it, and \cite{cakmak2020bayesian} studied Bayesian optimization of risk measures. Futhermore, \cite{shapiro2023bayesian} combined Bayesian estimation with distributional robustness to provide additional protection against misspecification of the parametric family.

For digital twins, DRO and BRO should consequently be viewed as complementary rather than mutually exclusive. DRO is compelling when certified protection against a clearly defined neighborhood of distributions is paramount. BRO is especially natural when observations arrive sequentially and the decision maker wants flexibility in choosing a risk attitude ranging from average-case to tail-focused or worst-case behavior. This separation between learning the model through the posterior and choosing an attitude toward posterior uncertainty is central to adaptive risk management.

\subsection{Distributionally robust MDPs and reinforcement learning}
Section~4 introduced the nominal Bellman recursion used for control and reinforcement learning. Here the sequential decision-making framework remains the same, but uncertainty about the transition model enters each Bellman backup. This makes robustness more delicate because adverse transition choices can compound over time. The classical robust-MDP formulations of \cite{iyengar2005robust} and \cite{nilim2005robust} optimize against transition kernels in uncertainty sets and show that, under rectangularity conditions, robust value functions still follow the dynamic-programming principle. Suppose that the one-stage cost $c(s,a)$ is known, and let $\mathcal{P}_{s,a}$ be a state--action rectangular uncertainty set of possible next-state distributions. In a discounted cost-minimization problem, the corresponding classical robust-MDP Bellman equation is
\begin{equation}
V^{\mathrm{R}}(s)
= \min_{a\in\mathcal{A}(s)}
\left\{
c(s,a) + \gamma
\sup_{p\in\mathcal{P}_{s,a}}
\mathbb{E}_{S'\sim p}
\big[V^{\mathrm{R}}(S')\big]
\right\}.
\label{eq:rmdp-bellman}
\end{equation}
Equation~\eqref{eq:rmdp-bellman} imposes uncertainty directly on transition distributions: the adversary selects $p$ from $\mathcal{P}_{s,a}$. Rectangularity permits this worst-case transition-kernel choice to be made locally inside each Bellman backup, which preserves dynamic programming. The same decoupling, however, can permit combinations of local worst cases that are statistically implausible and overly conservative.

A distributionally robust MDP (DR-MDP) adds a second layer of uncertainty by treating model parameters as random while leaving their probability law ambiguous. \cite{xu2012distributionally} optimize expected performance under the worst admissible law over uncertain reward and transition parameters and show that, under suitable structure, the problem reduces to a tractable robust MDP. Some later work also uses the term distributionally robust for direct ambiguity over transition distributions, but the distinction here is whether uncertainty is placed directly on transition kernels or on probability laws over uncertain models.

Distributionally robust reinforcement learning (DRRL) addresses the additional problem that the nominal transition model is learned from data. Much of the theory, however, has been developed in data regimes that are batch-like rather than continuously adaptive: \cite{panaganti2022sample} assume generative-model access, \cite{liu2022distributionally} use repeated simulator samples to construct a distributionally robust Q-learning estimator, \cite{shi2024distributionally} study a model-based offline-RL setting with a fixed historical dataset, and \cite{wang2023finite} established a finite-sample complexity bound for distributionally robust Q-learning with simulator access. Thus, a substantial part of the literature focuses on known-model planning, generative models, or offline data, with some exceptions for online RL. \cite{derman2020bayesian} use posterior uncertainty to adapt a robust Bellman equation and uncertainty set to new observations, and \cite{liang2024single} develop a fully model-free, multi-timescale algorithm that learns from a single incrementally observed trajectory.  \cite{lu2024interactive} studied an online-offline setting where data are collected interactively while learner interacts with the training environment only.

Even with these advances, the standard DR-MDP or DRRL objective can be conservative because it repeatedly optimizes against the least favorable distribution in a prespecified ambiguity set. New data may improve the nominal model and render a smaller ambiguity set, but neither the ambiguity-set update nor the desired level of pessimism is generated automatically by the data process; both require an additional calibration or adaptation rule. Moreover, the predominance of DRRL analyses based on planning, generative models, and offline data also leaves a gap for digital twins, where observations arrive sequentially and decisions affect what is learned next. Bayesian risk-aware MDPs (BR-MDPs) and reinforcement learning address these shortcomings by separating posterior learning from risk attitude and by allowing both to evolve during operation.

\subsection{Bayesian risk-aware MDPs and reinforcement learning}
Extending Bayesian risk optimization from a single-stage problem to an MDP is not immediate. A risk measure inserted into each Bellman backup must preserve a suitable dynamic structure, distinguish parameter uncertainty from ordinary transition randomness, and account for information that will become available at later stages. Let $\mu$ denote the current posterior over $\theta$ that parametrizes the transition kernel $P_\theta(\cdot|s,a)$, and let $\mathcal{B}(\mu;s,a,S')$ denote its Bayesian update after observing the transition to $S'$. A finite-horizon belief-state BR-MDP has the recursion
\begin{equation}
V_t^{\mathrm{BR}}(s,\mu)
= \min_{a\in\mathcal{A}(s)}
\rho_{\theta\sim\mu}\!\Bigg(
\mathbb{E}_{S'\sim P_{\theta}(\cdot\mid s,a)}\Big[
c_{\theta}(s,a,S') 
+V_{t+1}^{\mathrm{BR}}\!\left(
S',\mathcal{B}(\mu;s,a,S')
\right)
\Big]\Bigg).
\label{eq:brmdp-bellman}
\end{equation}
The inner expectation in Equation~\eqref{eq:brmdp-bellman} treats the transition of state-action conditional on $\theta$, whereas $\rho_{\theta\sim\mu}$ evaluates epistemic uncertainty under the posterior. Because the posterior is carried forward by $\mathcal{B}$, the augmented state $(s,\mu)$ makes the resulting policy adaptive to future learning part of the decision problem.

\cite{lin2022bayesian} introduced the finite-horizon belief-state BR-MDP in Equation~\eqref{eq:brmdp-bellman}. Posterior expectation recovers a Bayesian-average MDP, tail risk measures such as conditional value-at-risk yield intermediate protection, and limiting choices approach worst-case behavior. For conditional value-at-risk, they developed an efficient approximation based on an analytical value-function approximation and convexity.

A another modeling choice is to update the posterior only between episodes. In the episodic Bayesian optimal-control approach of \cite{shapiro2025episodic}, the posterior $\mu^{(k)}$ is formed at the beginning of episode $k$ and held fixed while a Bayesian-average control problem is solved and its policy is implemented. Observations collected during that episode are used to construct $\mu^{(k+1)}$ for the next episode. Consequently, the physical state remains part of the within-episode Bellman recursion, but the posterior is a fixed episode-level parameter rather than a state variable: the controller can react to changes in the physical state, but it does not adapt the policy to newly collected information until the beginning of the next episode. This formulation is computationally simpler and natural when a digital twin is synchronized periodically. The belief-state BR-MDP is yields a more adaptive policy because it values and reacts to information arriving in the future, but augmenting the state by a posterior increases computational complexity.

Relative to a distributionally robust MDP, these Bayesian formulations separate the statistical model and the evidence supplied by data from the decision criterion. In the belief-state BR-MDP, a risk functional additionally specifies how posterior uncertainty affects each Bellman backup; in the episodic Bayesian-average formulation, expectation is used and learning occurs at episode boundaries. As data accumulate and the posterior contracts, epistemic uncertainty decreases; in a BR-MDP, the gap among risk-neutral, tail-risk, and worst-case evaluations also narrows without an ad-hoc schedule for shrinking an ambiguity radius. A BR-MDP is therefore a natural step up for digital twins when within-horizon adaptation and flexible risk attitude are important.

\cite{wang2023bayesianrl} extended this perspective to infinite-horizon reinforcement learning with externally arriving observations. Their Bayesian risk-averse Q-learning method updates a posterior over transition models and converges to the true-environment optimum as all state--action transitions receive sufficient physical observations. 

Bayesian risk-aware RL also makes the interaction between robustness and exploration explicit. Online actions determine not only rewards but also future information. A pessimistic robust or distributionally robust policy may avoid poorly known state--action pairs precisely because their transitions are uncertain, even when visiting them would be informative; an optimistic policy may instead explore aggressively while exposing the physical system to poorly understood outcomes. Because posterior risk is represented separately from posterior learning, BR-MDP and Bayesian risk-aware RL can adjust this balance as evidence accumulates. \cite{song2026evolving} study this issue through quantile BR-MDPs. The posterior quantile level determines how model uncertainty enters the Bellman backup: lower-tail quantiles induce pessimism and robustness, whereas upper-tail quantiles induce optimism and encourage exploration. Because the posterior itself contracts, the effect of either choice weakens as evidence accumulates. They propose an adaptive quantile schedule that emphasizes robustness early and gradually increases exploration of less-visited state--action pairs, with sublinear Bayesian-regret guarantees relative to true-environment and robust Bayesian-risk benchmarks.


\section{Challenges and Future Directions}

Despite the recent advances reviewed above, significant challenges remain in applying digital twins to real-world decision-making problems. Here, we highlight several of these challenges and outline promising directions for future research. 

\subsection{Online optimization under nonstationarity and evolving simulators}
Existing streaming-data methods provide important first steps, but many assume that the physical system is stationary and that several simulation or optimization iterations can be completed between model updates. In practice, the physical system $M_t$ may itself evolve because of degradation, seasonal effects, behavioral responses, or changes in the operating environment. Consequently, discrepancies between $M_t$ and the digital-twin model $\widehat{M}_t$ may reflect both statistical estimation error and genuine distributional or structural change. An online decision method must therefore detect nonstationarity, discount or reweight obsolete information, update the twin, and track a moving optimum rather than repeatedly solve a sequence of unrelated static problems. Important questions include how quickly such methods can adapt after gradual drift or abrupt changes and how decision quality depends jointly on detection delay, model-update frequency, and optimization error. \par Multi-fidelity simulation offers a promising way to meet the associated computational demands. A digital twin may combine inexpensive reduced-order models, learned surrogates, and coarse simulations with a more accurate but costly high-fidelity simulator. Online decision making must then allocate computation across both candidate decisions and fidelity levels: low-fidelity models can support rapid screening and frequent policy updates, while high-fidelity evaluations can correct bias, validate promising decisions, and recalibrate the lower-fidelity models. 
Developing methods that jointly assimilate streaming physical data, select simulation fidelities, and track evolving optimal decisions is an important direction for digital-twin research.

\subsection{Endogenous data, active learning, and safe exploration}
Streaming data are not necessarily exogenous. In a digital twin, implemented decisions change the physical system, determine which states are visited, and influence which observations become available for the next model update. The resulting data can therefore be selectively sampled and statistically dependent on earlier decisions. If this feedback is ignored, the twin may become systematically accurate in frequently visited operating regimes while remaining poorly calibrated precisely where alternative decisions would take the system. \par This endogeneity links several methodological traditions. The Bayesian SGD method reviewed above addresses decision-dependent input distributions, experimental design and active learning choose observations for their value in reducing uncertainty, dual control recognizes that actions both regulate a system and reveal its dynamics, and online reinforcement learning balances exploitation with exploration. For digital twins, these viewpoints should be integrated with simulation: virtual experiments can screen informative or risky actions, whereas carefully chosen physical interventions supply evidence that simulation alone cannot provide. 

An important research direction is to develop safe policies that balance immediate operating performance, the cost of physical experimentation, and the value of better future decisions. Such policies should coordinate  physical data collection and multi-fidelity simulation; respect safety and operational constraints during exploration; and account for selection bias when updating the twin model $\widehat{M}_t$ from adaptively collected data. The resulting guarantees should characterize not only statistical learning rates or regret, but also how quickly the digital twin becomes trustworthy in decision-relevant regions of the state--action space.


A promising research direction is a multilayer framework in which simulation optimization handles operational design decisions, reinforcement learning adapts policies through experience, and control provides short-term feedback and stability. The main challenge is not merely connecting the algorithms, but ensuring that uncertainty and information are propagated consistently across the layers.

\subsection{Risk-aware autonomy}
As digital twins become more autonomous, erroneous decisions may be implemented before a human can intervene. Risk-aware methods must therefore do more than select a conservative policy: they should make uncertainty, safety margins, and the operational consequences of risk preferences transparent, and they should adapt protection as data accumulate. Bayesian risk formulations offer one mechanism because posterior learning and risk attitude are represented separately, but comparable adaptive calibration is needed for robust and distributionally robust approaches. \par Autonomy also magnifies the consequences of model misspecification and distribution shift. A twin may be confident for the wrong reasons when its model class omits an important failure mode, when sensors degrade, or when the deployed policy drives the system outside the region supported by data. Future methods should combine detection of model inadequacy with fallback policies, safety filters, and explicit triggers for reverting to a conservative mode or requesting human intervention. 

Evaluation standards should consequently extend beyond average reward or nominal optimality. Relevant guarantees might include constraint-violation probabilities, tail performance, recovery after detected model failure, and the reliability of uncertainty quantification under repeated model updates. Developing computationally scalable methods that provide such guarantees in real time is essential if risk-aware digital twins are to support high-stake autonomous decisions.

\subsection{Human in the loop}\par Even as digital twins become increasingly autonomous, human expertise remains essential for specifying objectives and constraints, recognizing novel operating regimes, and judging trade-offs that are difficult to encode in a reward or cost function. A human-in-the-loop digital twin should therefore support bidirectional learning: the twin communicates predictions, uncertainty, and the consequences of candidate actions in an interpretable form, while operators provide preferences, corrections, demonstrations, or overrides that update the decision model. These inputs should be distinguished from observations of the physical system so that model learning, preference learning, and uncertainty quantification are not conflated. Important research questions include when the twin should act autonomously, when it should defer to or query a person, how human feedback can be assimilated without double counting or introducing systematic bias, and how safety and performance guarantees can be maintained as both the model and human preferences evolve. Addressing these questions will require integrating simulation optimization and sequential decision making with preference learning and shared autonomy.

\section{Conclusions}
This review discusses decision-making challenges arising from digital twins. Unlike conventional simulation models, digital twins synchronize with physical systems through streaming data and repeated decisions. Therefore, data-driven optimization, sequential decision making, and risk-aware methodologies play important roles in applying digital twins to real-world problems.

The central challenge is that decisions must be made while the model is being learned and updated. In data-driven optimization, this leads to problems involving streaming input data, input uncertainty, and simultaneous allocation of data-collection and simulation resources. In sequential decision making, it leads to reinforcement learning and optimal control methods that update models and policies over time. In risk-aware decision making, it requires a treatment of model uncertainty that adapts as data accumulate.

Research in these areas is still developing. Digital twins are likely to create substantial demand for new methodologies that combine simulation, optimization, control, and learning. A unified decision-centric perspective may help connect these research communities and support the next generation of adaptive and risk-aware digital twins.

\section{Acknowledgement}
The author gratefully acknowledges the support by the Air Force Office of Scientific Research (AFOSR) under Grant FA9550-25-1-0310, the
National Science Foundation under Award ECCS-2419562, 
 and the NSF AI
Institute for Advances in Optimization under Grant NSF-2112533. The author used ChatGPT for polishing the writing and generating the figure in the paper.

\bibliographystyle{plain}
\bibliography{references}

@inproceedings{tirinzoni2019transfer,
  author    = {Tirinzoni, Andrea and Salvini, Mattia and Restelli, Marcello},
  title     = {Transfer of Samples in Policy Search via Multiple Importance Sampling},
  booktitle = {Proceedings of the 36th International Conference on Machine Learning},
  pages     = {6264--6274},
  year      = {2019},
  volume    = {97},
  series    = {Proceedings of Machine Learning Research},
  publisher = {PMLR}
}

@inproceedings{chandak2020optimizing,
  author    = {Chandak, Yash and Theocharous, Georgios and Shankar, Shiv and White, Martha and Mahadevan, Sridhar and Thomas, Philip},
  title     = {Optimizing for the Future in Non-Stationary {MDP}s},
  booktitle = {Proceedings of the 37th International Conference on Machine Learning},
  pages     = {1414--1425},
  year      = {2020},
  volume    = {119},
  series    = {Proceedings of Machine Learning Research},
  publisher = {PMLR}
}

@inproceedings{liu2023reusing,
  author    = {Liu, Vincent and Chandak, Yash and Thomas, Philip and White, Martha},
  title     = {Asymptotically Unbiased Off-Policy Policy Evaluation when Reusing Old Data in Nonstationary Environments},
  booktitle = {Proceedings of the 26th International Conference on Artificial Intelligence and Statistics},
  pages     = {5474--5492},
  year      = {2023},
  volume    = {206},
  series    = {Proceedings of Machine Learning Research},
  publisher = {PMLR}
}

@inproceedings{queeney2021generalized,
  author    = {Queeney, James and Paschalidis, Ioannis Ch. and Cassandras, Christos G.},
  title     = {Generalized Proximal Policy Optimization with Sample Reuse},
  booktitle = {Advances in Neural Information Processing Systems},
  volume    = {34},
  pages     = {11909--11919},
  year      = {2021},
  publisher = {Curran Associates, Inc.}
}

@inproceedings{sutton1990integrated,
  author    = {Sutton, Richard S.},
  title     = {Integrated Architectures for Learning, Planning, and Reacting Based on Approximating Dynamic Programming},
  booktitle = {Proceedings of the Seventh International Conference on Machine Learning},
  pages     = {216--224},
  year      = {1990},
}

@inproceedings{janner2019trust,
  author    = {Janner, Michael and Fu, Justin and Zhang, Marvin and Levine, Sergey},
  title     = {When to Trust Your Model: Model-Based Policy Optimization},
  booktitle = {Advances in Neural Information Processing Systems},
  volume    = {32},
  pages     = {12498--12509},
  year      = {2019},
  publisher = {Curran Associates, Inc.}
}

@inproceedings{peng2018simtoreal,
  author    = {Peng, Xue Bin and Andrychowicz, Marcin and Zaremba, Wojciech and Abbeel, Pieter},
  title     = {Sim-to-Real Transfer of Robotic Control with Dynamics Randomization},
  booktitle = {Proceedings of the 2018 IEEE International Conference on Robotics and Automation},
  pages     = {3803--3810},
  year      = {2018},
  publisher = {IEEE},
}

@inproceedings{chebotar2019closing,
  author    = {Chebotar, Yevgen and Handa, Ankur and Makoviychuk, Viktor and Macklin, Miles and Issac, Jan and Ratliff, Nathan and Fox, Dieter},
  title     = {Closing the Sim-to-Real Loop: Adapting Simulation Randomization with Real World Experience},
  booktitle = {Proceedings of the 2019 IEEE International Conference on Robotics and Automation},
  pages     = {8973--8979},
  year      = {2019},
  publisher = {IEEE},
}

@article{voogd2023reinforcement,
  author  = {Voogd, Kevin L. and Allamaa, Jean Pierre and Alonso-Mora, Javier and Son, Tong Duy},
  title   = {Reinforcement Learning from Simulation to Real World Autonomous Driving Using Digital Twin},
  journal = {IFAC-PapersOnLine},
  volume  = {56},
  number  = {2},
  pages   = {1510--1515},
  year    = {2023},
}

@inproceedings{deisenroth2011pilco,
  author    = {Deisenroth, Marc Peter and Rasmussen, Carl Edward},
  title     = {{PILCO}: A Model-Based and Data-Efficient Approach to Policy Search},
  booktitle = {Proceedings of the 28th International Conference on Machine Learning},
  pages     = {465--472},
  year      = {2011}
}

@article{hewing2020learning,
  author  = {Hewing, Lukas and Wabersich, Kim P. and Menner, Marcel and Zeilinger, Melanie N.},
  title   = {Learning-Based Model Predictive Control: Toward Safe Learning in Control},
  journal = {Annual Review of Control, Robotics, and Autonomous Systems},
  volume  = {3},
  pages   = {269--296},
  year    = {2020},
}

@inproceedings{koller2018safe,
  author    = {Koller, Torsten and Berkenkamp, Felix and Turchetta, Matteo and Krause, Andreas},
  title     = {Learning-Based Model Predictive Control for Safe Exploration},
  booktitle = {Proceedings of the 57th IEEE Conference on Decision and Control},
  pages     = {6059--6066},
  year      = {2018},
  publisher = {IEEE},
}

@inproceedings{wabersich2020bayesian,
  author    = {Wabersich, Kim Peter and Zeilinger, Melanie N.},
  title     = {Bayesian Model Predictive Control: Efficient Model Exploration and Regret Bounds Using Posterior Sampling},
  booktitle = {Proceedings of the 2nd Conference on Learning for Dynamics and Control},
  pages     = {455--464},
  year      = {2020},
  volume    = {120},
  series    = {Proceedings of Machine Learning Research},
  publisher = {PMLR}
}

@article{he2024stochastic,
  author    = {He, Linyun and Shanbhag, Uday V. and Song, Eunhye},
  title     = {Stochastic Approximation for Multi-period Simulation Optimization with Streaming Input Data},
  journal   = {ACM Transactions on Modeling and Computer Simulation},
  volume    = {34},
  number    = {2},
  articleno = {6},
  numpages  = {27},
  year      = {2024},
  month     = apr,
}

@article{nasem2024digital,
  author      = {{National Academies of Sciences, Engineering, and Medicine}},
  title       = {Foundational Research Gaps and Future Directions for Digital Twins},
  journal = {National Academies Press},
  year        = {2024},
  address     = {Washington, DC}
}

@article{kritzinger2018digital,
  author  = {Kritzinger, Werner and Karner, Matthias and Traar, Georg and Henjes, Jan and Sihn, Wilfried},
  title   = {Digital Twin in Manufacturing: A Categorical Literature Review and Classification},
  journal = {IFAC-PapersOnLine},
  volume  = {51},
  number  = {11},
  pages   = {1016--1022},
  year    = {2018}
}

@article{jones2020characterising,
  author  = {Jones, David and Snider, Chris and Nassehi, Aydin and Yon, Jason and Hicks, Ben},
  title   = {Characterising the Digital Twin: A Systematic Literature Review},
  journal = {CIRP Journal of Manufacturing Science and Technology},
  volume  = {29},
  pages   = {36--52},
  year    = {2020}
}

@inproceedings{heluany2023survey,
  author    = {Heluany, Jessica B. and Gkioulos, Vasileios},
  title     = {Survey on Digital Twins: From Concepts to Applications},
  booktitle = {Proceedings of the 18th International Conference on Availability, Reliability and Security},
  year      = {2023},
}

@article{wu2023comprehensive,
  author  = {Wu, Honghai and Ji, Pengwei and Ma, Huahong and Xing, Ling},
  title   = {A Comprehensive Review of Digital Twin from the Perspective of Total Process: Data, Models, Networks and Applications},
  journal = {Sensors},
  volume  = {23},
  number  = {19},
  pages   = {8306},
  year    = {2023},
}

@inproceedings{biller2022simulation,
  author    = {Biller, Bahar and Jiang, Xiaoning and Yi, Jinglei and Venditti, Paula and Biller, Stephan},
  title     = {Simulation: The Critical Technology in Digital Twin Development},
  booktitle = {Proceedings of the 2022 Winter Simulation Conference},
  pages     = {1340--1355},
  year      = {2022}
}

@inproceedings{matta2023digital,
  author    = {Matta, Andrea and Lugaresi, Giacomo},
  title     = {Digital Twins: Features, Models, and Services},
  booktitle = {Proceedings of the 2023 Winter Simulation Conference},
  pages     = {46--60},
  year      = {2023}
}

@inproceedings{zhou2024tutorial,
  author = {Zhou, Enlu},
  title = {Data-Driven Simulation Optimization in the Age of Digital Twins: Challenges and Developments},
  booktitle = {Proceedings of the 2024 Winter Simulation Conference},
  pages = {31--45},
  year = {2024}
}

@inproceedings{li2020o2des,
author    = {Li, Haobin and Cao, Xiaofeng and Sharma, P. and Lee, Loo Hay and Chew, Ek Peng},
  title     = {Framework of O2DES.NET Digital Twins for Next Generation Ports and Warehouse Solutions},
  booktitle = {Proceedings of the 2020 Winter Simulation Conference},
  pages     = {3188--3199},
  year      = {2020}
}

@book{fu2015handbook,
  editor    = {Fu, Michael C.},
  title     = {Handbook of Simulation Optimization},
  publisher = {Springer},
  series    = {International Series in Operations Research \& Management Science},
  volume    = {216},
  year      = {2015}
}

@article{chick2001input,
  author  = {Chick, Stephen E.},
  title   = {Input Distribution Selection for Simulation Experiments: Accounting for Input Uncertainty},
  journal = {Operations Research},
  volume  = {49},
  number  = {5},
  pages   = {744--758},
  year    = {2001}
}

@article{zouaoui2003accounting,
  author  = {Zouaoui, Faker and Wilson, James R.},
  title   = {Accounting for Parameter Uncertainty in Simulation Input Modeling},
  journal = {IIE Transactions},
  volume  = {35},
  pages   = {781--792},
  year    = {2003}
}

@article{zouaoui2004accounting,
  author  = {Zouaoui, Faker and Wilson, James R.},
  title   = {Accounting for Input-Model and Input-Parameter Uncertainties in Simulation},
  journal = {IIE Transactions},
  volume  = {36},
  number  = {11},
  pages   = {1135--1151},
  year    = {2004}
}

@article{xie2014bayesian,
  author  = {Xie, Weiwei and Nelson, Barry L. and Barton, Russell R.},
  title   = {A Bayesian Framework for Quantifying Uncertainty in Stochastic Simulation},
  journal = {Operations Research},
  volume  = {62},
  number  = {6},
  pages   = {1439--1452},
  year    = {2014}
}

@article{barton2014quantifying,
  author  = {Barton, Russell R. and Nelson, Barry L. and Xie, Weiwei},
  title   = {Quantifying Input Uncertainty via Simulation Confidence Intervals},
  journal = {INFORMS Journal on Computing},
  volume  = {26},
  number  = {1},
  pages   = {74--87},
  year    = {2014}
}

@article{gao2017robust,
  author  = {Gao, Siyang and Xiao, Haibo and Zhou, Enlu and Chen, Wei},
  title   = {Robust Ranking and Selection with Optimal Computing Budget Allocation},
  journal = {Automatica},
  volume  = {81},
  pages   = {30--36},
  year    = {2017}
}

@article{xiao2018simulation,
  author  = {Xiao, Haibo and Gao, Siyang},
  title   = {Simulation Budget Allocation for Selecting the Top-$m$ Designs with Input Uncertainty},
  journal = {IEEE Transactions on Automatic Control},
  volume  = {63},
  number  = {9},
  pages   = {3127--3134},
  year    = {2018}
}

@article{xiao2020optimal,
  author  = {Xiao, Haibo and Gao, Siyang and Lee, Loo Hay},
  title   = {Optimal Computing Budget Allocation for Complete Ranking with Input Uncertainty},
  journal = {IISE Transactions},
  volume  = {52},
  number  = {5},
  pages   = {489--499},
  year    = {2020}
}

@article{fan2020distributionally,
  author  = {Fan, Wei and Hong, L. Jeff and Zhang, Xiaowei},
  title   = {Distributionally Robust Selection of the Best},
  journal = {Management Science},
  volume  = {66},
  number  = {1},
  pages   = {190--208},
  year    = {2020}
}

@article{wu2018bayesian,
  author  = {Wu, Di and Zhu, Helin and Zhou, Enlu},
  title   = {A Bayesian Risk Approach to Data-Driven Stochastic Optimization: Formulations and Asymptotics},
  journal = {SIAM Journal on Optimization},
  volume  = {28},
  number  = {2},
  pages   = {1588--1612},
  year    = {2018}
}

@article{cakmak2021solving,
  author  = {Cakmak, Sait and Wu, Di and Zhou, Enlu},
  title   = {Solving {B}ayesian Risk Optimization via Nested Stochastic Gradient Estimation},
  journal = {IISE Transactions},
  pages   = {1081--1093},
  year    = {2021}
}

@inproceedings{cakmak2020bayesian,
  author    = {Cakmak, Sait and Astudillo, Raul and Frazier, Peter I. and Zhou, Enlu},
  title     = {Bayesian Optimization of Risk Measures},
  booktitle = {Advances in Neural Information Processing Systems},
  year      = {2020}
}

@article{shapiro2023bayesian,
  author  = {Shapiro, Alexander and Zhou, Enlu and Lin, Yifan},
  title   = {Bayesian Distributionally Robust Optimization},
  journal = {SIAM Journal on Optimization},
  volume  = {33},
  number  = {2},
  pages   = {1279--1304},
  year    = {2023}
}

@article{wu2024datadriven,
  author  = {Wu, Di and Wang, Yuhao and Zhou, Enlu},
  title   = {Data-Driven Ranking and Selection under Input Uncertainty},
  journal = {Operations Research},
  volume  = {72},
  number  = {2},
  pages   = {781--795},
  year    = {2024}
}

@inproceedings{wang2022fixed,
  author    = {Wang, Yuhao and Zhou, Enlu},
  title     = {Fixed Budget Ranking and Selection with Streaming Input Data},
  booktitle = {Proceedings of the 2022 Winter Simulation Conference},
  pages     = {3027--3038},
  year      = {2022}
}

@article{liu2024bayesian,
  author  = {Liu, Tianyi and Lin, Yifan and Zhou, Enlu},
  title   = {Bayesian Stochastic Gradient Descent for Stochastic Optimization with Streaming Input Data},
  journal = {SIAM Journal on Optimization},
  volume  = {34},
  number  = {1},
  pages   = {389--418},
  year    = {2024}
}

@inproceedings{wu2017ranking,
  author    = {Wu, Di and Zhou, Enlu},
  title     = {Ranking and Selection under Input Uncertainty: A Budget Allocation Formulation},
  booktitle = {Proceedings of the 2017 Winter Simulation Conference},
  pages     = {2245--2256},
  year      = {2017}
}

@inproceedings{xu2020joint,
  author    = {Xu, Jie and Zheng, Zeyu and Glynn, Peter W.},
  title     = {Joint Resource Allocation for Input Data Collection and Simulation},
  booktitle = {Proceedings of the 2020 Winter Simulation Conference},
  pages     = {2126--2137},
  year      = {2020}
}

@inproceedings{kim2022optimizing,
  author    = {Kim, Taeho and Song, Eunhye},
  title     = {Optimizing Input Data Acquisition for Ranking and Selection: A View Through the Most Probable Best},
  booktitle = {Proceedings of the 2022 Winter Simulation Conference},
  pages     = {2258--2269},
  year      = {2022}
}

@article{wang2025ranking,
  author = {Wang, Yuhao and Zhou, Enlu},
  title = {Data-Driven Ranking and Selection With Simultaneous Input Data Collection and Simulation},
  journal = {Naval Research Logistics},
  volume = {72},
  number = {8},
  pages = {1170--1183},
  year = {2025}
}

@article{shen2021ranking,
  author  = {Shen, Haihui and Hong, L. Jeff and Zhang, Xiaowei},
  title   = {Ranking and Selection with Covariates for Personalized Decision Making},
  journal = {INFORMS Journal on Computing},
  volume  = {33},
  number  = {4},
  pages   = {1500--1519},
  year    = {2021}
}

@inproceedings{li2020context,
  author    = {Li, Haobin and Lam, Henry and Peng, Yijie},
  title     = {Context-Dependent Ranking and Selection under a Bayesian Framework},
  booktitle = {Proceedings of the 2020 Winter Simulation Conference},
  pages     = {2060--2070},
  year      = {2020}
}

@article{li2022efficient,
  author = {Li, Haidong and Lam, Henry and Peng, Yijie},
  title = {Efficient Learning for Clustering and Optimizing Context-Dependent Designs},
  journal = {Operations Research},
  volume = {72},
  number = {2},
  pages = {617--638},
  year = {2024}
}

@article{cakmak2024contextual,
  author = {Cakmak, Sait and Wang, Yuhao and Gao, Siyang and Zhou, Enlu},
  title = {Contextual Ranking and Selection with Gaussian Processes and Optimal Computing Budget Allocation},
  journal = {ACM Transactions on Modeling and Computer Simulation},
  volume = {34},
  number = {2},
  articleno = {8},
  numpages = {24},
  year = {2024}
}

@article{hong2019offline,
  author  = {Hong, L. Jeff and Jiang, Guangxin},
  title   = {Offline Simulation Online Application: A New Framework of Simulation-Based Decision Making},
  journal = {Asia-Pacific Journal of Operational Research},
  volume  = {36},
  pages   = {1940015},
  year    = {2019}
}

@article{jiang2020online,
  author  = {Jiang, Guangxin and Hong, L. Jeff and Nelson, Barry L.},
  title   = {Online Risk Monitoring Using Offline Simulation},
  journal = {INFORMS Journal on Computing},
  volume  = {32},
  number  = {2},
  pages   = {356--375},
  year    = {2020}
}

@book{sutton2018reinforcement,
  author    = {Sutton, Richard S. and Barto, Andrew G.},
  title     = {Reinforcement Learning: An Introduction},
  edition   = {2},
  publisher = {MIT Press},
  year      = {2018}
}

@book{bertsekas2017dynamic,
  author    = {Bertsekas, Dimitri P.},
  title     = {Dynamic Programming and Optimal Control},
  publisher = {Athena Scientific},
  year      = {2017}
}

@misc{schena2024reinforcement,
  author = {Schena, Lorenzo and Marques, Pedro and Poletti, Romain and Ahizi, Samuel and Van den Berghe, Jan and Mendez, Miguel A.},
  title  = {Reinforcement Twinning: From Digital Twins to Model-Based Reinforcement Learning},
  year   = {2024},
  note   = {arXiv:2311.03628}
}

@inproceedings{wang2023bayesianrl,
  author    = {Wang, Yuhao and Zhou, Enlu},
  title     = {Bayesian Risk-Averse {Q}-Learning with Streaming Observations},
  booktitle = {Advances in Neural Information Processing Systems},
  year      = {2023},
}

@article{shapiro2025episodic,
  author = {Shapiro, Alexander and Zhou, Enlu and Lin, Yifan and Wang, Yuhao},
  title = {Episodic {B}ayesian Optimal Control with Unknown Randomness Distributions},
  journal = {Operations Research},
  volume = {74},
  number = {3},
  pages = {1437--1455},
  year = {2025}
}

@article{lin2025reusing,
author = {Lin, Yifan and Wang, Yuhao and Zhou, Enlu},
  title = {Reusing Historical Trajectories in Natural Policy Gradient via Importance Sampling: Convergence and Convergence Rate},
  journal = {Operations Research},
  volume = {73},
  number = {6},
  pages = {3010--3026},
  year = {2025}
}

@article{iyengar2005robust,
  author  = {Iyengar, Garud N.},
  title   = {Robust Dynamic Programming},
  journal = {Mathematics of Operations Research},
  volume  = {30},
  number  = {2},
  pages   = {257--280},
  year    = {2005}
}

@article{nilim2005robust,
  author  = {Nilim, Arnab and El Ghaoui, Laurent},
  title   = {Robust Control of Markov Decision Processes with Uncertain Transition Matrices},
  journal = {Operations Research},
  volume  = {53},
  number  = {5},
  pages   = {780--798},
  year    = {2005}
}

@misc{song2026evolving,
  author = {Song, Meichen and Wang, Yuhao and Zhou, Enlu},
  title  = {Evolving Robustness--Exploration Trade-off in Online Reinforcement Learning via Quantile {B}ayesian Risk {MDP}s},
  year   = {2026},
  note   = {arXiv:2605.24345}
}

@inproceedings{lin2022bayesian,
  author    = {Lin, Yifan and Ren, Yuxuan and Zhou, Enlu},
  title     = {Bayesian Risk Markov Decision Processes},
  booktitle = {Advances in Neural Information Processing Systems},
  volume    = {35},
  year      = {2022},
  url       = {https://proceedings.neurips.cc/paper_files/paper/2022/hash/6f7d90b1198fec96defd80b5ebd5bc81-Abstract-Conference.html}
}

@article{delage2010distributionally,
  author  = {Delage, Erick and Ye, Yinyu},
  title   = {Distributionally Robust Optimization Under Moment Uncertainty with Application to Data-Driven Problems},
  journal = {Operations Research},
  volume  = {58},
  number  = {3},
  pages   = {595--612},
  year    = {2010},
}

@article{wiesemann2014distributionally,
  author  = {Wiesemann, Wolfram and Kuhn, Daniel and Sim, Melvyn},
  title   = {Distributionally Robust Convex Optimization},
  journal = {Operations Research},
  volume  = {62},
  number  = {6},
  pages   = {1358--1376},
  year    = {2014},
}

@article{esfahani2018data,
  author  = {Mohajerin Esfahani, Peyman and Kuhn, Daniel},
  title   = {Data-Driven Distributionally Robust Optimization Using the {Wasserstein} Metric: Performance Guarantees and Tractable Reformulations},
  journal = {Mathematical Programming},
  volume  = {171},
  pages   = {115--166},
  year    = {2018},
}

@article{xu2012distributionally,
  author  = {Xu, Huan and Mannor, Shie},
  title   = {Distributionally Robust Markov Decision Processes},
  journal = {Mathematics of Operations Research},
  volume  = {37},
  number  = {2},
  pages   = {288--300},
  year    = {2012},
}

@inproceedings{derman2020bayesian,
  author    = {Derman, Esther and Mankowitz, Daniel and Mann, Timothy and Mannor, Shie},
  title     = {A Bayesian Approach to Robust Reinforcement Learning},
  booktitle = {Proceedings of the 35th Conference on Uncertainty in Artificial Intelligence},
  series    = {Proceedings of Machine Learning Research},
  volume    = {115},
  pages     = {648--658},
  year      = {2020},
  publisher = {PMLR}
}

@inproceedings{liu2022distributionally,
  author    = {Liu, Zijian and Bai, Qinxun and Blanchet, Jose and Dong, Perry and Xu, Wei and Zhou, Zhengqing and Zhou, Zhengyuan},
  title     = {Distributionally Robust {Q}-Learning},
  booktitle = {Proceedings of the 39th International Conference on Machine Learning},
  series    = {Proceedings of Machine Learning Research},
  volume    = {162},
  pages     = {13623--13643},
  year      = {2022},
  publisher = {PMLR}
}

@inproceedings{panaganti2022sample,
  author    = {Panaganti, Kishan and Kalathil, Dileep},
  title     = {Sample Complexity of Robust Reinforcement Learning with a Generative Model},
  booktitle = {Proceedings of the 25th International Conference on Artificial Intelligence and Statistics},
  series    = {Proceedings of Machine Learning Research},
  volume    = {151},
  pages     = {9582--9602},
  year      = {2022},
  publisher = {PMLR}
}

@inproceedings{liang2024single,
  author    = {Liang, Zhipeng and Ma, Xiaoteng and Blanchet, Jose and Yang, Jun and Zhang, Jiheng and Zhou, Zhengyuan},
  title     = {Single-Trajectory Distributionally Robust Reinforcement Learning},
  booktitle = {Proceedings of the 41st International Conference on Machine Learning},
  series    = {Proceedings of Machine Learning Research},
  volume    = {235},
  pages     = {29644--29666},
  year      = {2024},
  publisher = {PMLR}
}

@article{blanchet2019quantifying,
  author  = {Blanchet, Jose and Murthy, Karthyek},
  title   = {Quantifying Distributional Model Risk via Optimal Transport},
  journal = {Mathematics of Operations Research},
  volume  = {44},
  number  = {2},
  pages   = {565--600},
  year    = {2019},
}

@article{gao2023distributionally,
  author  = {Gao, Rui and Kleywegt, Anton J.},
  title   = {Distributionally Robust Stochastic Optimization with {Wasserstein} Distance},
  journal = {Mathematics of Operations Research},
  volume  = {48},
  number  = {2},
  pages   = {603--655},
  year    = {2023},
}

@article{shi2024distributionally,
  author  = {Shi, Laixi and Chi, Yuejie},
  title   = {Distributionally Robust Model-Based Offline Reinforcement Learning with Near-Optimal Sample Complexity},
  journal = {Journal of Machine Learning Research},
  volume  = {25},
  number  = {200},
  pages   = {1--91},
  year    = {2024}
}

@article{wang2024optimal,
  author  = {Wang, Yuhao and Zhou, Enlu},
  title   = {Optimal Computing Budget Allocation for Data-Driven Ranking and Selection},
  journal = {INFORMS Journal on Optimization},
  volume  = {7},
  number  = {1},
  pages   = {1--19},
  year    = {2024},
}

@unpublished{li2026bayesianmpc,
  author = {Li, Yingke and Lin, Yifan and Zhou, Enlu and Zhang, Fumin},
  title = {Bayesian Risk-Averse Model Predictive Control with Consistency and Stability Guarantees},
  note = {Submitted; arXiv:2511.21871v1},
  year = {2025}
}

@inproceedings{wang2023finite,
  author = {Wang, Shengbo and Si, Nian and Blanchet, Jose and Zhou, Zhengyuan},
  title = {A Finite Sample Complexity Bound for Distributionally Robust Q-Learning},
  booktitle = {Proceedings of the 26th International Conference on Artificial Intelligence and Statistics},
  series = {Proceedings of Machine Learning Research},
  volume = {206},
  pages = {3370--3398},
  year = {2023},
  publisher = {PMLR}
}

@inproceedings{lu2024interactive,
  author = {Lu, Miao and Zhong, Han and Zhang, Tong and Blanchet, Jose},
  title = {Distributionally Robust Reinforcement Learning with Interactive Data Collection: Fundamental Hardness and Near-Optimal Algorithms},
  booktitle = {Advances in Neural Information Processing Systems},
  volume = {37},
  pages = {12528--12580},
  year = {2024}
}

@article{hong2021review,
  author = {Hong, L. Jeff and Fan, Weiwei and Luo, Jun},
  title = {Review on Ranking and Selection: A New Perspective},
  journal = {Frontiers of Engineering Management},
  volume = {8},
  number = {3},
  pages = {321--343},
  year = {2021},
}

@article{du2024contextual,
  author  = {Du, Jianzhong and Gao, Siyang and Chen, Chun-Hung},
  title   = {A Contextual Ranking and Selection Method for Personalized Medicine},
  journal = {Manufacturing and Service Operations Management},
  volume  = {26},
  number  = {1},
  pages   = {167--181},
  year    = {2024},
}

@article{keslin2024ranking,
  author = {Keslin, Gregory and Nelson, Barry L. and Pagnoncelli, Bernardo and Plumlee, Matthew and Rahimian, Hamed},
  title = {Ranking and Contextual Selection},
  journal = {Operations Research},
  volume = {73},
  number = {5},
  pages = {2695--2707},
  year = {2025}
}

@article{zhang2025topm,
  author  = {Zhang, Gongbo and Chen, Sihua and Huang, Kuihua and Peng, Yijie},
  title   = {Efficient Learning for Selecting Top-$m$ Context-Dependent Designs},
  journal = {IEEE Transactions on Automation Science and Engineering},
  volume  = {22},
  pages   = {3210--3225},
  year    = {2025},
}

@article{li2026quadratic,
  author  = {Li, Dongyang and Chew, Ek Peng and Li, Haobin and Y{\"u}cesan, Enver and Chen, Chun-Hung},
  title   = {Efficient Simulation Budget Allocation for Contextual Ranking and Selection with Quadratic Models},
  journal = {European Journal of Operational Research},
  volume  = {328},
  number  = {3},
  pages   = {862--876},
  year    = {2026},
}

@article{sadana2025survey,
  author  = {Sadana, Utsav and Chenreddy, Abhilash and Delage, Erick and Forel, Alexandre and Frejinger, Emma and Vidal, Thibaut},
  title   = {A Survey of Contextual Optimization Methods for Decision-Making under Uncertainty},
  journal = {European Journal of Operational Research},
  volume  = {320},
  number  = {2},
  pages   = {271--289},
  year    = {2025},
}

@article{lam2022subsampling,
  author = {Lam, Henry and Qian, Huajie},
  title = {Subsampling to Enhance Efficiency in Input Uncertainty Quantification},
  journal = {Operations Research},
  volume = {70},
  number = {3},
  pages = {1891--1913},
  year = {2022}
}

@article{zhu2020risk,
  author = {Zhu, Helin and Liu, Tianyi and Zhou, Enlu},
  title = {Risk Quantification in Stochastic Simulation under Input Uncertainty},
  journal = {ACM Transactions on Modeling and Computer Simulation},
  volume = {30},
  number = {1},
  pages = {1--24},
  year = {2020}
}

@inproceedings{liu2021bayesian,
  author = {Liu, Tianyi and Lin, Yifan and Zhou, Enlu},
  title = {A {B}ayesian Approach to Online Simulation Optimization with Streaming Input Data},
  booktitle = {Proceedings of the 2021 Winter Simulation Conference},
  pages = {1--12},
  year = {2021}
}
\end{document}